%% file: paper.tex
\input 11layout
\input macro

\centerline{\biggbf The Interleaved Multichromatic Number of a Graph}

\bigskip\bigskip
\centerline{\it Valmir C. Barbosa}

\medskip
\centerline{Universidade Federal do Rio de Janeiro}
\centerline{Programa de Engenharia de Sistemas e Computa\c c\~ao, COPPE}
\centerline{Caixa Postal 68511}
\centerline{21941-972 Rio de Janeiro - RJ, Brazil}
\centerline{\tt valmir@cos.ufrj.br}

\bigskip\bigskip
\centerline{\bf Abstract}

\medskip
\noindent
For $k\ge 1$, we consider interleaved $k$-tuple colorings of the nodes of a
graph, that is, assignments of $k$ distinct natural numbers to each node
in such a way that nodes that are connected by an edge receive numbers that
are strictly alternating between them with respect to the relation $<$. If
it takes at least $\chi_{\rm int}^k(G)$ distinct numbers to provide graph
$G$ with such a coloring, then the interleaved multichromatic number of $G$
is $\chi_{\rm int}^*(G)=\inf_{k\ge 1}\chi_{\rm int}^k(G)/k$ and is known to
be given by a function of the simple cycles of $G$ under acyclic orientations
if $G$ is connected.\footnote{$^1$}{V. C. Barbosa and E. Gafni, {\it ACM
Trans.\ on Programming Languages and Systems\/} {\bf 11} (1989), 562--584.}
This paper contains a new proof of this result. Unlike the original proof, the
new proof makes no assumptions on the connectedness of $G$, nor does it resort
to the possible applications of interleaved $k$-tuple colorings and their
properties.

\bigskip\bigskip
\noindent
{\bf Keywords:} Chromatic number, multichromatic number, interleaved
multichromatic number, acyclic orientations.

\vfill\eject
\bigbeginsection 1. Introduction

Let $G=(N,E)$ be an undirected graph with
$\vert E\vert\ge 1$. For $k\ge 1$, a {\it $k$-tuple coloring\/} of $G$ is an
assignment of $k$ distinct natural numbers ({\it colors\/}) to each of the nodes
of $G$ in such a way that no two nodes connected by an edge are assigned an
identical color. The least number of colors with which $G$ can be $k$-tuple
colored is its {\it $k$-chromatic number}, denoted by $\chi^k(G)$. For $k=1$,
the $1$-chromatic number of $G$ is the graph's {\it chromatic number}, denoted
by $\chi(G)$.

There is also a notion of ``efficiency'' that goes along $k$-tuple colorings of
$G$, which is the notion that, for two distinct positive integers $k_1$ and
$k_2$, it is more ``efficient'' to $k_1$-tuple color $G$ than to $k_2$-tuple
color it if and only if $\chi^{k_1}(G)/k_1<\chi^{k_2}(G)/k_2$. This notion is
formalized by the definition of the {\it multichromatic number\/} of $G$ [5],
denoted by $\chi^*(G)$ and given by
$$\chi^*(G)=\inf_{k\ge 1}{{\chi^k(G)}\over{k}}.\eqno(1)$$
Because the infimum in (1) can be shown to be always attained [2], a $k$ for
which the ratio in (1) is minimum gives a most ``efficient'' $k$-tuple
coloring of $G$.

One special case of $k$-tuple colorings that has great practical interest is
the case of {\it interleaved $k$-tuple colorings}, defined as follows. For
$i\in N$, let $c_i^1,\ldots,c_i^k$ be the $k$ colors assigned to node $i$.
This $k$-tuple coloring is interleaved if and only if either
$c_i^1<c_j^1<\cdots<c_i^k<c_j^k$ or $c_j^1<c_i^1<\cdots<c_j^k<c_i^k$ for all
$(i,j)\in E$. Every $1$-tuple coloring is interleaved, and the counterparts for
interleaved colorings of the graph's $k$-chromatic and multichromatic numbers
are, respectively, its {\it interleaved $k$-chromatic number\/} and
{\it interleaved multichromatic number}. These are denoted, respectively,
by $\chi_{\rm int}^k(G)$ and $\chi_{\rm int}^*(G)$. Following (1), the latter
is defined as
$$\chi_{\rm int}^*(G)=\inf_{k\ge 1}{{\chi_{\rm int}^k(G)}\over{k}}.\eqno(2)$$
Clearly, $\chi^*(G)\le\chi_{\rm int}^*(G)\le\chi(G)$.

For connected $G$, the infimum in (2) has also been shown to be always attained
and given as follows [1]. Let $\Kappa$ be the set of all simple cycles of $G$.
For $\kappa\in\Kappa$, let $\kappa^+$ and $\kappa^-$ be the two possible
traversal directions of $\kappa$. If $\omega$ is an {\it acyclic orientation\/}
of $G$ (that is, an assignment of directions to edges that forms no directed
cycles), then we denote by $m(\kappa^+,\omega)$ the number of edges in $\kappa$
oriented by $\omega$ in the $\kappa^+$ direction. Likewise for
$m(\kappa^-,\omega)$. What has been shown for nonempty $\Kappa$ is that
$${{1}\over{\chi_{\rm int}^*(G)}}=
\max_{\omega\in\Omega(G)}\min_{\kappa\in\Kappa}
{{\min\bigl\{m(\kappa^+,\omega),m(\kappa^-,\omega)\bigr\}}
\over{\vert\kappa\vert}},\eqno(3)$$
where $\Omega(G)$ is the set of all acyclic orientations of $G$ and
$\vert\kappa\vert$ is the size of $\kappa$, that is,
$\vert\kappa\vert=m(\kappa^+,\omega)+m(\kappa^-,\omega)$. If $\Kappa=\emptyset$
(i.e., $G$ is a tree), then $\chi_{\rm int}^*(G)=2$.

In [1], $1/\chi_{\rm int}^*(G)$ is analyzed as a measure of concurrency (or
parallelism) for certain distributed computations over systems that can be
represented by connected graphs, and (3) is obtained through proofs that make
heavy use of both the graph's connectedness and properties of those distributed
computations. What we do in this paper is to provide an alternative proof for
(3), one that does not make assumptions on the connectedness of $G$ and is
built from ``first principles,'' that is, independently of the possible
applications of the concept of ``efficient'' interleaved $k$-tuple colorings.
This is done in Section 2, and in Section 3 we give concluding remarks and
argue that (3) provides both ``efficient'' interleaved $k$-tuple colorings and
optimal $1$-tuple colorings with a common underlying relationship to the set
$\Omega(G)$ of the graph's acyclic orientations.

\bigbeginsection 2. A new proof

Our proof is based on the {\it lexicographic product\/} [4] of $G$ and the
complete graph on $k$ nodes. This lexicographic product is another undirected
graph, denoted by $G^k$ and having $k$ nodes for each node in $G$. For $i\in N$,
we denote the nodes of $G^k$ by $i^1,\ldots,i^k$. The edges of $G^k$ are
deployed in such a way that, for fixed $i$, an edge exists connecting every
pair of nonidentical nodes of $i^1,\ldots,i^k$, and for $(i,j)\in E$, an edge
exists between each one of $i^1,\ldots,i^k$ and each one of $j^1,\ldots,j^k$.
Clearly, every $k$-tuple coloring of $G$ corresponds to a $1$-tuple coloring
of $G^k$, and conversely.

Now consider the set $\Omega(G^k)$ of the acyclic orientations of $G^k$. We
say that $\varphi\in\Omega(G^k)$ is a {\it layered\/} orientation of $G^k$ if
and only if every edge $(i^{k_1},j^{k_2})$ for which $k_1>k_2$ is oriented
from $i^{k_1}$ to $j^{k_2}$, and furthermore every edge $(i^{k_1},j^{k_1})$
gets the same orientation for $k_1=1,\ldots,k$. If we visualize $G^k$ as
$k$ instances of $G$ stacked upon one another in such a way that the instance
at the bottom contains superscript-$1$ nodes, the one right above it
superscript-$2$ nodes, and so on, then what a layered acyclic orientation of
$G^k$ does is to provide the $k$ instances of $G$ with identical acyclic
orientations from $\Omega(G)$ and orient all other edges (the ones that connect
nodes from different instances of $G$) from top to bottom. The set of all
layered acyclic orientations of $G^k$ is denoted by $\Omega^L(G^k)$.

In $G^k$, we refer to each of the aforementioned instances of $G$ as a
{\it layer}, and number each layer with the superscript of the nodes that it
contains. Given a layered acyclic orientation, every maximal directed path
starts in layer $k$ (because that is where all {\it sources\/} are, these being
nodes whose adjacent edges are all oriented outward) and ends in layer $1$
(which contains all {\it sinks}, that is, nodes whose adjacent edges are all
oriented inward).

For a generic undirected graph $H$ and for an acyclic orientation $\omega$ of
$H$, let $l_\omega$ be the number of nodes in the longest directed path in $H$
according to $\omega$. Also, let a $1$-tuple coloring of $H$ be called
{\it monotonic\/} with respect to an acyclic orientation $\omega$ of $H$ if and
only if the colors it assigns decrease along any directed path according to
$\omega$. If $\Omega(H)$ is the set of acyclic orientations of $H$, then we
have the following, for whose proof the reader is referred to [3].

\proclaim Lemma 1. To every $\omega\in\Omega(H)$ there corresponds a monotonic
$1$-tuple coloring of $H$ with at most $l_\omega$ colors. Conversely, to every
$1$-tuple coloring of $H$ by a total of $c$ colors there corresponds an
orientation $\omega\in\Omega(H)$ with respect to which the coloring is monotonic
and such that $l_\omega\le c$.

By Lemma 1, the chromatic number of graph $H$ is given by the number of nodes
on the longest directed path that is shortest among all possible acyclic
orientations of $H$ [3]. In the case of $G^k$, and considering that
$\chi^k(G)=\chi(G^k)$, we have
$\chi^k(G)=\min_{\varphi\in\Omega(G^k)}l_\varphi$. Lemma 1 is the basis
of our first supporting result, which states that it suffices to look at the
orientations in $\Omega^L(G^k)\subseteq\Omega(G^k)$ to obtain
$\chi_{\rm int}^k(G)$.

\proclaim Lemma 2. For $k\ge 1$,
$\chi_{\rm int}^k(G)=\min_{\varphi\in\Omega^L(G^k)}l_\varphi$.

\proof
Let $\varphi$ be a layered acyclic orientation of $G^k$. By Lemma 1, $G^k$ can
be colored monotonically by a $1$-tuple coloring that employs a total of at most
$l_\varphi$ colors. If the instance of $(i,j)\in E$ within each layer is
oriented by $\varphi$ from the corresponding instance of $i$ to the
corresponding instance of $j$, then by definition of a layered orientation a
directed path exists from $i^k$ to $j^1$ that alternates instances of $i$ and
$j$ along its way (that is, a path like
$i^k\to j^k\to i^{k-1}\to j^{k-1}\to\cdots\to i^1\to j^1$). It then follows from
the monotonicity of $1$-tuple colorings of $G^k$ that such a coloring
corresponds to an interleaved $k$-tuple coloring of $G$, thus
$$\chi_{\rm int}^k(G)\le\min_{\varphi\in\Omega^L(G^k)}l_\varphi.$$

Conversely, consider an interleaved $k$-tuple coloring of $G$ by a total of $c$
colors. This coloring corresponds to possibly several $1$-tuple colorings of
$G^k$ by the same $c$ colors. Of these, and for $i\in N$, consider the $1$-tuple
coloring that assigns the highest of the $k$ colors of $i$ to $i^k$, then the
next highest to $i^{k-1}$, and so on. By Lemma 1, an orientation $\varphi$ of
$G^k$ exists with respect to which this $1$-tuple coloring of $G^k$ is monotonic
and as such orients all edges of $G^k$ from the node with higher color to the
one with lower. For an arbitrary edge $(i,j)$ of $G$, this, together with the
fact that the assumed $k$-tuple coloring of $G$ is interleaved, implies that
$\varphi$ is layered. In addition, Lemma 1 also ensures that $l_\varphi\le c$,
hence
$$\chi_{\rm int}^k(G)\ge\min_{\varphi\in\Omega^L(G^k)}l_\varphi.$$
\endproclaim

Our next supporting result is a statement on the morphology of longest directed
paths according to layered acyclic orientations of $G^k$. What it states is that
such a directed path's every edge is either confined to a layer or joins two
adjacent layers. In the latter case, it also states that, if such an edge exists
from $i^{k_1}$ to $j^{k_1-1}$, then a directed edge also exists from $j^{k_1}$
to $i^{k_1}$.

\proclaim Lemma 3. For $k\ge k_1\ge k_2\ge 1$ and $\varphi\in\Omega^L(G^k)$,
let $p$ be a longest directed path according to $\varphi$, and let $i^{k_1}$
and $j^{k_2}$ be such that a directed edge exists on $p$ from $i^{k_1}$ to
$j^{k_2}$. Then, for $\ell=1,\ldots,k$, either $k_1=k_2$ and a directed edge
exists from $i^\ell$ to $j^\ell$, or $k_1=k_2+1$ with $i\neq j$ and a directed
edge exists from $j^\ell$ to $i^\ell$.

\proof
Note, first, that $k_1\le k_2+1$, since $k_1>k_2+1$ contradicts the hypothesis
that $p$ is longest. If $k_1=k_2$, then the lemma follows directly from the
definition of a layered acyclic orientation of $G^k$. We then assume that
$k_1=k_2+1$.

What remains to be shown is that two specific scenarios cannot happen if
$\varphi$ is layered. The first one is that $i=j$, that is, $p$ contains an edge
directed from $i^{k_1}$ to $i^{k_1-1}$. That this cannot happen follows from the
fact that replacing this edge on $p$ with either the directed path
$i^{k_1}\to r^{k_1}\to i^{k_1-1}$ or the directed path
$i^{k_1}\to r^{k_1-1}\to i^{k_1-1}$ for some $(i,r)\in E$ yields a directed path
longer than $p$. Because $\vert E\vert\ge 1$, at least one of the two
possibilities is certain to occur, so we have a contradiction.

The second scenario corresponds to $(i,j)\in E$ when every layer $\ell$ contains
an edge directed from $i^\ell$ to $j^\ell$. In this case, an edge exists
directed from $i^{k_1}$ to $j^{k_1-1}$ on $p$. Once again, this edge can be
replaced with the directed path $i^{k_1}\to j^{k_1}\to i^{k_1-1}\to j^{k_1-1}$
to yield a directed path longer than $p$, also a contradiction.
\endproclaim

By Lemma 3, a longest directed path $p$ in $G^k$ according to a layered acyclic
orientation $\varphi$ looks like the following. It starts at a source in layer
$k$, winds its way through the $k$ layers, and ends at a sink in layer $1$. In
addition, for each of the $k-1$ edges that lead from one layer to the next, say
for example an edge from $i^{k_1}$ to $j^{k_1-1}$ for $1<k_1\le k$, an edge
exists in layer $\ell$ directed from $j^\ell$ to $i^\ell$ for $\ell=1,\ldots,k$.

If we let $\omega\in\Omega(G)$ be the acyclic orientation of $G$ that
corresponds to the orientation of a layer by $\varphi$ (the same for all
layers), then $p$ can be seen to leave a {\it trace\/} on $G$ oriented by
$\omega$ as well. This trace is an undirected path in $G$ (not necessarily a
simple one) that can be traversed by following $\omega$ on the edges that
correspond to $p$ being confined to a layer, and by going against $\omega$ on
the $k-1$ edges that correspond to the transition of $p$ from one layer to the
next. Of course, leaving such a trace on $G$ oriented by $\omega$ is no
prerogative of longest directed paths, but rather applies to any other directed
path whose edges are either confined to a same layer or join adjacent layers as
described in Lemma 3.

Next we give our main result, in which (3) is established. Recall that $\Kappa$
stands for the set of all simple cycles of $G$.

\proclaim Theorem 4. If $G$ is a forest, then $\chi_{\rm int}^*(G)=2$.
Otherwise,
$$\chi_{\rm int}^*(G)=
\min_{\omega\in\Omega(G)}\max_{\kappa\in\Kappa}
{{\vert\kappa\vert}
\over{\min\bigl\{m(\kappa^+,\omega),m(\kappa^-,\omega)\bigr\}}}.$$

\proof
By definition of $\chi_{\rm int}^*(G)$ (cf.\ (2)) and by Lemma 2, we have
$$\eqalign{
\chi_{\rm int}^*(G)&=\inf_{k\ge 1}{{\chi_{\rm int}^k(G)}\over{k}}\cr
&=\inf_{k\ge 1}{{1}\over{k}}\min_{\varphi\in\Omega^L(G^k)}l_\varphi.\cr
}$$
Also, the definition of layered acyclic orientations of $G^k$ implies that
there exists a one-to-one correspondence between such orientations and the
acyclic orientations of $G$. If, for $k\ge 1$ and $\omega\in\Omega(G)$, we let
$\varphi_\omega^k$ be the corresponding orientation in $\Omega^L(G^k)$, then we
can further write
$$\eqalign{
\chi_{\rm int}^*(G)
&=\inf_{k\ge 1}{{1}\over{k}}\min_{\omega\in\Omega(G)}l_{\varphi_\omega^k}\cr
&=\min_{\omega\in\Omega(G)}\inf_{k\ge 1}{{l_{\varphi_\omega^k}}\over{k}}.\cr
}$$

If $G$ is a forest, then the theorem follows trivially (this is also a
consequence of the facts that $\chi^*(G)=\chi(G)=2$ if $G$ is a forest and that
$\chi^*(G)\le\chi_{\rm int}^*(G)\le\chi(G)$). If not, then it suffices to argue
that
$$\max_{\kappa\in\Kappa}
{{\vert\kappa\vert}
\over{\min\bigl\{m(\kappa^+,\omega),m(\kappa^-,\omega)\bigr\}}}
=\inf_{k\ge 1}{{l_{\varphi_\omega^k}}\over{k}}$$
for all $\omega\in\Omega(G)$.

Let $\hat\Kappa\subseteq\Kappa$ be the set of all simple cycles $\kappa$ of $G$
for which
$\vert\kappa\vert/\min\bigl\{m(\kappa^+,\omega),m(\kappa^-,\omega)\bigr\}$ is
maximum over $\Kappa$, and for all $\kappa\in\hat\Kappa$ assume that
$m(\kappa^+,\omega)\ge m(\kappa^-,\omega)$.
For $k\ge 1$, every $\kappa\in\hat\Kappa$ can be seen to give rise to directed
paths in $G^k$ whose edges either lead to other nodes in the same layer or else
descend to the next layer as described in Lemma 3. As we discussed following the
presentation of that lemma, each such path leaves a trace on $G$ oriented by
$\omega$. One possibility is the path that we denote by $p_k(\kappa,\omega)$ and
whose trace has the following characteristics. It starts at the origin of the
longest segment of $\kappa$ oriented by $\omega$ in the $\kappa^+$ direction,
then winds its way around $\kappa$ a number of times given by
$t=\bigl\lfloor(k-1)/m(\kappa^-,\omega)\bigr\rfloor$, and finally goes around
$\kappa$ one last time until $(k-1)\bmod m(\kappa^-,\omega)$ additional edges
have been traversed against $\omega$ and one last segment has been traversed in
agreement with $\omega$. Along this last traversal around $\kappa$, let
$\epsilon^+$ denote the number of edges traversed in agreement with $\omega$ and
$\epsilon^-$ of those traversed against $\omega$. Evidently,
$\epsilon^-\le m(\kappa^-,\omega)-1$.

Note that this last traversal around $\kappa$ starts, like $p_k(\kappa,\omega)$,
at the origin of the longest segment of $\kappa$ oriented by $\omega$ in the
$\kappa^+$ direction. If, in addition, it is chosen to be longest among all
possibilities, then we have
$$\epsilon^+\ge(1+\epsilon^-){{m(\kappa^+,\omega)}\over{m(\kappa^-,\omega)}}.$$
Equivalently,
$$\epsilon^+m(\kappa^-,\omega)\ge(1+\epsilon^-)m(\kappa^+,\omega),$$
or yet,
$$\displaylines{
\qquad\qquad
\bigl(1+tm(\kappa^+,\omega)+tm(\kappa^-,\omega)+\epsilon^++\epsilon^-\bigr)
m(\kappa^-,\omega)\hfill\cr
\hfill\ge\bigl(1+tm(\kappa^-,\omega)+\epsilon^-\bigr)
\bigl(m(\kappa^+,\omega)+m(\kappa^-,\omega)\bigr),\qquad\qquad\cr
}$$
which yields
$$\eqalign{
{{1+tm(\kappa^+,\omega)+tm(\kappa^-,\omega)+\epsilon^++\epsilon^-}
\over{1+tm(\kappa^-,\omega)+\epsilon^-}}
&\ge{{m(\kappa^+,\omega)+m(\kappa^-,\omega)}
\over{m(\kappa^-,\omega)}}\cr
&={{\vert\kappa\vert}\over{m(\kappa^-,\omega)}}.\cr
}$$

If we let $l_k(\kappa,\omega)$ denote the number of nodes of
$p_k(\kappa,\omega)$, then the left-hand side of the latter inequality is equal
to $l_k(\kappa,\omega)/k$. Considering further that
$l_{\varphi_\omega^k}\ge l_k(\kappa,\omega)$, we obtain
$${{l_{\varphi_\omega^k}}\over{k}}\ge
{{\vert\kappa\vert}\over{m(\kappa^-,\omega)}}.$$
Also, for sufficiently large $k$, a longest directed path in $G^k$ (comprising
$l_{\varphi_\omega^k}$ nodes, by definition) leaves a trace on $G$ oriented
by $\omega$ which, after a preamble that comprises a number of nodes that is
bounded as $k$ grows, coincides with one of the members, say $\hat\kappa$, of
$\hat\Kappa$. This is ensured by the fact that, for every
$m(\hat\kappa^-,\omega)$ times that the path has to change layers, the greatest
contribution to its number of nodes is $\vert\hat\kappa\vert$, by definition of
$\hat\Kappa$.

So we can write, by Lemma 3,
$${{l_{\varphi_\omega^k}}\over{k}}=
{{tm(\hat\kappa^+,\omega)+tm(\hat\kappa^-,\omega)+o(k)}
\over{tm(\hat\kappa^-,\omega)+o(k)}},$$
following the standard use of $o(k)$ to signify that
$\lim_{k\to\infty}o(k)/k=0$, for some $t$ that grows without bounds along with
$k$. Clearly, then,
$$\lim_{k\to\infty}{{l_{\varphi_\omega^k}}\over{k}}=
{{\vert\hat\kappa\vert}\over{m(\hat\kappa^-,\omega)}},$$
which concludes the proof.
\endproclaim

\bigbeginsection 3. Discussion

This paper contains a new proof of (3). This new proof, in contrast to the
original proof that appears in [1], does not make assumptions on the
connectedness of $G$, nor does it resort to the specifics of possible
applications of the concept of interleaved $k$-tuple colorings. Unlike the
proof of [1], however, the proof in this paper does not establish that the
infimum in (2) is always achieved, which is indeed the case. For this,
the original proof continues to be the source, as establishing that property
is achieved independently of the connectedness of $G$.

As expressed in Theorem 4, $\chi_{\rm int}^*(G)$ corresponds to the minimum,
taken over all of the acyclic orientations of $G$, of a function of certain
undirected paths of $G$ under those acyclic orientations. Specifically, if
for $\omega\in\Omega(G)$ we let $P(\omega)$ comprise every path $p$ in $G$
such that $p$ contains all nodes in a simple cycle and starts at a node that
is a source according to $\omega$ if only the edges on the simple cycle are
considered, then $\chi_{\rm int}^*(G)$ can be written as
$$\chi_{\rm int}^*(G)=\min_{\omega\in\Omega(G)}\max_{p\in P(\omega)}
{{l(p)}\over{1+c(p,\omega)}},\eqno(4)$$
where $l(p)$ is the number of nodes of $p$ and $c(p,\omega)$ is the number
of edges on $p$ that are oriented by $\omega$ contrary to the traversal of
$p$ from its starting node.

Interestingly, the expression in (4) is the same as the aforementioned
consequence of Lemma 1 that
$$\chi(G)=\min_{\omega\in\Omega(G)}\max_{p\in P(\omega)}l(p),\eqno(5)$$
if only we now let $P(\omega)$ stand for the set of directed paths of $G$
according to $\omega$, since in this case $c(p,\omega)=0$ for all
$\omega\in\Omega(G)$ and all $p\in P(\omega)$.

Since by (4) and (5), respectively, it is possible to assess a graph's
interleaved multichromatic number and chromatic number by looking at functions
of the graph's paths as oriented by the acyclic orientations of $G$, one
interesting question is whether $\chi^*(G)$, the multichromatic number of $G$,
is also amenable to a similar characterization. Unlike the other two cases,
however, there does not appear to be a result like Lemma 2, which establishes
the fundamental one-to-one correspondence between the acyclic orientations of
$G^k$ that matter and those of $G$. Finding such a characterization for
$\chi^*(G)$ remains then an open question.

\beginsection  Acknowledgments

The author acknowledges partial support from CNPq, CAPES, the PRONEX initiative
of Brazil's MCT under contract 41.96.0857.00, and a FAPERJ BBP grant.

\bigbeginsection References

{\frenchspacing

\medskip
\item{1.} V. C. Barbosa and E. Gafni,
``Concurrency in heavily loaded neighborhood-constrained systems,''
{\it ACM Trans. on Programming Languages and Systems\/} {\bf 11} (1989),
562--584.

\medskip
\item{2.} F. H. Clarke and R. E. Jamison,
``Multicolorings, measures and games on graphs,''
{\it Discrete Mathematics\/} {\bf 14} (1976),
241--245.

\medskip
\item{3.} R. W. Deming,
``Acyclic orientations of a graph and chromatic and independence numbers,''
{\it J. of Combinatorial Theory B\/} {\bf 26} (1979),
101--110.

\medskip
\item{4.} D. Geller and S. Stahl,
``The chromatic number and other functions of the lexicographic product,''
{\it J. of Combinatorial Theory B\/} {\bf 19} (1975),
87--95.

\medskip
\item{5.} S. Stahl,
``$n$-tuple colorings and associated graphs,''
{\it J. of Combinatorial Theory B\/} {\bf 20} (1976),
185--203.

}

\bye

%% file: 11layout.tex
\magnification=1100
\hsize=6.5 true in
\vsize=9.0 true in
\overfullrule=0pt
\baselineskip=15pt
\tolerance=10000
\hbadness=200

%% file: macro.tex
\font\tensc=cmcsc10
\newfam\scfam
\textfont\scfam=\tensc
\def\sc{\fam\scfam\tensc}
 
\font\bigbf=cmbx10 scaled\magstep1

\font\biggbf=cmbx10 scaled\magstep2

\font\bigggbf=cmbx10 scaled\magstep3

\font\bigit=cmti10 scaled\magstep1

\def\Kappa{{\rm K}}

\def\sectionbreak{
\bigskip\vskip\parskip}

\def\bigsectionbreak{
\bigskip\bigskip\vskip\parskip}

\outer\def\beginsection#1\par
{\sectionbreak
\message{#1}\leftline{\bf#1}\nobreak\smallskip\noindent}

\outer\def\bigbeginsection#1\par
{\sectionbreak
\message{#1}\leftline{\bigbf#1}\nobreak\medskip\noindent}

\def\currentsection{\firstmark}

\outer\def\biggbeginsection#1\par
{\bigsectionbreak
\message{#1}\leftline{\biggbf#1}
\mark{#1}\nobreak\bigskip\noindent}

\outer\def\bigitbeginsection#1\par
{\sectionbreak
\message{#1}\leftline{\bigit#1}\nobreak\medskip\noindent}

\outer\def\longbigbeginsection#1 #2\par#3\par
{\sectionbreak
\message{#1 #2 #3}
\halign{##\hfil&##\hfil\cr
{\bigbf#1\ }&{\bigbf#2}\cr
&{\bigbf#3}\cr
}\nobreak\medskip\noindent}

\outer\def\longbiggbeginsection#1 #2\par#3\par
{\bigsectionbreak
\message{#1 #2 #3}
\halign{##\hfil&##\hfil\cr
{\biggbf#1\ }&{\biggbf#2}\cr
&{\biggbf#3}\cr
}\mark{#1 #2 #3}\nobreak\bigskip\noindent}

\outer\def\longbiggbeginappendix#1 #2 #3\par#4\par
{\bigsectionbreak
\message{#1 #2 #3 #4}
\halign{##\hfil&##\hfil\cr
{\biggbf#1\ #2\ }&{\biggbf#3}\cr
&{\biggbf#4}\cr
}\mark{#1 #2 #3 #4}\nobreak\bigskip\noindent}

\def\rightheadline{\hfil{\it\currentsection}\ \ \ \ \ \ {\rm \folio}}

\def\currentchapter{}

\def\leftheadline{{\rm \folio}\ \ \ \ \ \ {\it\currentchapter}\hfil}

\newcount\titlepageno

\def\setheadline{\headline=
{\ifnum\titlepageno=\pageno{\hfil}
\else{\ifodd\pageno{\rightheadline}\else{\leftheadline}\fi}
\fi}}

\def\skipifeven
{\ifodd\pageno{}\else\advancepageno\fi}

\outer\def\beginchapter#1. #2\par
{\vfill\eject\skipifeven\titlepageno=\pageno
\def\currentchapter{Cap\'\i tulo #1. #2}
\topinsert\vskip 0.25\vsize\endinsert
\hrule\medskip
\rightline{\bigggbf Cap\'\i tulo #1}
\medskip
\rightline{\bigggbf#2}
\medskip\hrule\bigskip\bigskip}

\outer\def\shortbeginchapter#1\par
{\vfill\eject\skipifeven\titlepageno=\pageno
\def\currentchapter{#1}
\mark{#1}
\topinsert\vskip 0.25\vsize\endinsert
\hrule\medskip
\rightline{\bigggbf#1}
\medskip\hrule\bigskip\bigskip}

\def\proof
{\noindent{\bf Proof: }}

\def\endproclaim
{{\unskip\nobreak\hfil\penalty50\hskip3.33pt\hbox{}\nobreak\hfil
\vrule height5pt width5pt depth0pt
\parfillskip=0pt \finalhyphendemerits=0 \medskip}}

\def\LaTeX{{\rm L\kern-.36em\raise.3ex\hbox{\sc a}\kern-.15em%
    T\kern-.1667em\lower.7ex\hbox{E}\kern-.125emX}}

\def\mytilde{\kern-.5pt\lower3pt\hbox{\char'176}\kern.5pt}